\documentclass[runningheads]{llncs}
\usepackage{graphicx}
\usepackage{latexsym}
\usepackage{amssymb}
\usepackage{amsmath}
\usepackage{mathrsfs}
\usepackage{enumerate}
\usepackage{booktabs}
\usepackage{bm}
\usepackage{xytree}

\DeclareFontEncoding{LS1}{}{}
\DeclareFontSubstitution{LS1}{stix2}{m}{n}

\newcommand{\blackdiamond}{{\text{\usefont{LS1}{stix2bb}{m}{it}\symbol{"E0}}}}

\begin{document}
\title{An Intermediate Logic Contained in Medvedev's Logic with Disjunction Property}
\titlerunning{Strong Union Logic}
\author{Zhicheng Chen}
\authorrunning{Z. Chen}
\institute{Peking University, Beijing, China}
\maketitle              
\begin{abstract}
Let \textbf{SU} be the superintuitionistic logic defined by the axiom
\(\boldsymbol{su} = ((\neg p\to q)\land(\neg q\to p) \rightarrow r \vee s) \to ( p \rightarrow r) \vee(q \rightarrow s)\), or equivalently, by Andrew's axiom. It is easy to  check that \textbf{SU} is contained in Medvedev's logic and contains both Kreisel-Putnam logic  and Scott logic.
We show that on \textbf{S4} frames, $\boldsymbol{su}$ corresponds to a certain first-order property, called the ``strong union'' property.
The strong completeness of \textbf{SU}, with respect to the class of \textbf{S4} frames enjoying this property, is proved. Furthermore, we demonstrate that \textbf{SU} has the disjunction property. As a result, \textbf{SU} stands as the strongest logic currently known below Medvedev's logic that has both an axiomatization and the disjunction property.
\end{abstract}

\section{Introduction}

Medvedev's logic \textbf{ML}, introduced by Y. Medvedev in 1962 \cite{Medvedev1962finiteproblems,Medvedev1966interpretation}, can be defined as the logic of the class of frames in the following form (using the intuitionistic propositional language):
\[\left\langle \wp^*(X), \supseteq\right\rangle\]
where $X$ is a non-empty finite set and $\wp^*(X)$ collects the non-empty subsets of $X$. 

\cite{Levin1969syntactic,Maksimova1986disjunction} proved that \textbf{ML} is a maximal intermediate logic with the ``disjunction property''. Additionally, \cite{Prucnal1979twoproblems} revealed that \textbf{ML} is ``structurally complete''. 
It is known that \textbf{ML} cannot be finitely axiomatized \cite{MaksimovaSS1979}, but whether \textbf{ML} is recursively axiomatizable (or equivalently, \textbf{ML} is decidable) remains a long-standing open problem.


Since \textbf{ML} contains intuitionistic propositional logic \textbf{IPL}, to approach the axiomatization of \textbf{ML}, one wants to find as strong as possible axioms of \textbf{ML} which go beyond \textbf{IPL}.

Among the known logics between \textbf{IPL} and \textbf{ML} in the literature, the most representative are the Kreisel-Putnam logic (\textbf{KP}) \cite{KreiselPutnam1957} and Scott logic (\textbf{SL}) \cite{KreiselPutnam1957}. They are obtained by adding respectively the following axioms to intuitionistic logic:
\begin{align*}
    \boldsymbol{kp} = (\neg p \rightarrow q \vee r) \to (\neg p \rightarrow q) \vee(\neg p \rightarrow r),\\
    \boldsymbol{sa} = ((\neg\neg p\to p)\to p\lor\neg p)\to\neg p\lor\neg\neg p.
\end{align*}

Both logics have the disjunction property \cite{KreiselPutnam1957} and are decidable \cite{Gabbay1970decidabilityKP,anderson_superconstructive_1972}. Moreover, \cite{minari_extension_1986} studied the logic \textbf{KPS} obtained by adding both $\boldsymbol{kp}$ and $\boldsymbol{sa}$ to \textbf{IPL}, proving that it has the disjunction property and is decidable.

Although it is known that there are axioms stronger than $\boldsymbol{kp}+\boldsymbol{sa}$ in \textbf{ML}, such as a class of axioms given in \cite{miglioli_results_1989} and the Andrew's axiom mentioned on page 107 of \cite{gabbay_semantical_2013}:
\[\boldsymbol{aa} = ((\neg p\to q) \rightarrow r \vee s) \to ((\neg p\to q) \rightarrow r) \vee(\neg\neg p \rightarrow s).\]
However, no literature has yet studied the properties of the logical systems obtained by adding these axioms to \textbf{IPL}.

In this paper, we will study the properties of the logic (denoted as \textbf{SU}) obtained by adding the following axiom to \textbf{IPL}:
\[\boldsymbol{su} = ((\neg p\to q)\land(\neg q\to p) \rightarrow r \vee s) \to ( p \rightarrow r) \vee(q \rightarrow s).\]

It is not difficult to see that $\boldsymbol{su}$ is an axiom of \textbf{ML}, and (on the basis of \textbf{IPL}) $\boldsymbol{su}$ is stronger than $\boldsymbol{aa}$. In fact, it is pointed out by Yifeng Ding that on top of \textbf{IPL}, $\boldsymbol{su}$ is equivalent to $\boldsymbol{aa}$. Nonetheless, our novel contribution is to show that \textbf{SU} is Kripke complete and has the disjunction property. As a result, \textbf{SU} stands as the currently strongest logic beneath \textbf{ML} that has both a finite axiomatization and the disjunction property.

We end this section by discussing the meaning of the $\boldsymbol{su}$ axiom.

In the Medvedev frame $\left\langle \wp^*(X), \supseteq\right\rangle$, any two points $w,v\in \wp^*(X)$ possess an infimum, which is their set union $w\cup v$. Moreover, it is easy to see that the following set-theoretic proposition (denoted as ($\star$)) holds: for any  $t\in \wp^*(w\cup v)$, $t$ belongs to $ \wp^*(w)$, or $t$ is in $\wp^*(v)$, or there are $w'\in \wp^*(w)$ and $v'\in \wp^*(v)$ with $t=w\cup v$.


Although it may be difficult to express ``any two points have an infimum'' in the intuitionistic language, combining ($\star$) it is not difficult to prove that Medvedev frames have the following property: for any  point $s$ in the frame, and any successors $w$ and $v$ of $s$, there exists  a successor $u$ of $s$ such that

\begin{center}
    u\text{ sees }w\text{ and }v\text{, and for any  successor }t\text{ of }u,\\
    t\text{ has a common successor with }w \text{ or } t\text{ has a common successor with }v \text{.}
\end{center}

This property is denoted as (Uni) (``Uni'' stands for union) in \cite{Holliday2017tolmf}. This property is actually expressible in the intuitionistic language—it corresponds precisely to the $\boldsymbol{kp}$ axiom on reflexive transitive frames.

However, there is more that can be expressed in the intuitionistic language regarding the infimum operation in Medvedev frames. Using ($\star$) it is not difficult to prove that Medvedev frames have the following property: for any  point $s$ in the frame, and any successors $w$ and $v$ of $s$, there exists  a successor $u$ of $s$ such that

\begin{center}
u\text{ sees }w\text{ and }v\text{, and for any  successor }t\text{ of }u, \\
t\text{ has a common successor with }w\text{ or }v\text{ sees }t\text{, and vice versa.}
\end{center}

We call this property the ``strong union'' property, denoted as (su). We will prove that this property (on reflexive transitive frames) corresponds precisely to the $\boldsymbol{su}$ axiom.

\section{Axiom $\boldsymbol{su}$ and its frame correspondence}

\begin{definition}
Let $PV=\{p_0,p_1,...\}$ be the countably infinite set of propositional variables. Let $Form$ collect all the formulas built up from propositional variables and $\perp$, using binary connectives $\land,\lor,\to$. 
 
Define the unary connective $\neg$: $\neg \alpha=(\alpha\to\,\perp).$ The order of precedence of the connectives is $\neg>\lor=\land>\,\to$. 
\end{definition}

\begin{definition}
$\boldsymbol{su} = ((\neg p\to q)\land(\neg q\to p) \rightarrow r \vee s) \to ( p \rightarrow r) \vee(q \rightarrow s)$.
\end{definition}

Below we show that $\boldsymbol{su}$ is equivalent to Andrew's axiom on the basis of $\textbf{IPL}$.
\begin{lemma}\label{su equals aa}
    On top of $\textbf{IPL}$, the following axioms are equivalent: 
    \begin{enumerate}
        \item $\boldsymbol{su}$;
        \item $\boldsymbol{aa^+} = ((p\to q) \rightarrow r \vee s) \to (( p\to q) \rightarrow r) \vee(\neg p \rightarrow s)$;
        \item $\boldsymbol{aa} = ((\neg p\to q) \rightarrow r \vee s) \to ((\neg p\to q) \rightarrow r) \vee(\neg\neg p \rightarrow s)$.
    \end{enumerate}
\end{lemma}
\begin{proof}
    It's easy to see that 1 is stronger than 2 and that 2 is stronger than 3. Now let $\,\vdash_{\textbf{A}}$ be the syntactic consequence relation obtained by adding $\boldsymbol{aa}$ (as an axiom schema) to the Hilbert system of $\textbf{IPL}$. We want to show $\,\vdash_{\textbf{A}}\boldsymbol{su}$. In fact, we have:
    
    $(\neg p\to q)\land(\neg q\to p) \rightarrow r \vee s$

    \noindent$\vdash_{\textbf{A}} (\neg p\to q)\to((\neg q\to p) \rightarrow r \vee s)$ \hfil\quad 

    \noindent$\vdash_{\textbf{A}} (\neg p\to q)\to ((\neg q\to p) \rightarrow r)) \vee (( \neg\neg q\rightarrow s)$ \hfill\quad ($\boldsymbol{aa}$)

    \noindent$\vdash_{\textbf{A}} (\neg\neg p\to((\neg q\to p) \rightarrow r)) \vee ((\neg p\to q)\to  (\neg\neg q\rightarrow s))$ \hfill\quad ($\boldsymbol{aa}$)

    \noindent$\vdash_{\textbf{A}} (p \rightarrow r) \vee (q\rightarrow s)$ \hfill\quad 
    
\qed
\end{proof}

\begin{definition}
For any Kripke frame \( \mathfrak{F} = \langle W, R \rangle \) and \( X, Y \subseteq W \), define:
\begin{itemize}
    \item[] $\blackdiamond_{\mathfrak{F}}(X)=\{w\in W\mid \textit{there is an } x \in X \textit{ with }xRw\}$,
    \item[] $\Diamond_{\mathfrak{F}}(X)=\{w\in W\mid \textit{there is an } x \in X \textit{ with }wRx\}$,
    \item[] $\Box_{\mathfrak{F}}(X) = \{ w \in W \mid$ for any $v$ with $wRv$, v is in $X \}$,
    \item[] $-_{\mathfrak{F}}(X) = W \setminus X$.
\end{itemize}
\end{definition}
(Note: for convenience, in the following we will omit the subscript “\(\mathfrak{F}\)” and parentheses.)

\begin{definition}[Strong union property]
~
\begin{enumerate}
\item For any  $n \in \omega^*$, any \textbf{S4} frame (i.e. reflexive and transitive frame) $\mathfrak{F} = \langle W, R \rangle$, any $z \in W$, and any $n$-tuple $x_0,\ldots,x_{n-1}$ in $W$, we say that $z$ strongly unites (/is the strong union of) $x_0,\ldots,x_{n-1}$ in $\mathfrak{F}$, if and only if
\[\text{for any  } i \in n\text{, } z R x_i \text{ and } z \in \Box(\blackdiamond\{x_i\} \cup \bigcup_{i' \in n \setminus \{i\}}\Diamond \blackdiamond\{x_{i'}\})\text{.}\]

In particular, for any  $z,x_0 \in W$, $z$ strongly unites $x_0$ iff $z R x_0$ and $x_0 R z$.

\item For any  \textbf{S4} frame $\mathfrak{F} = \langle W, R \rangle$, and any $n \in \omega^*$, define:
\begin{itemize}
    \item $\mathfrak{F}$ satisfies $(su_n) \Leftrightarrow \text{for any } w \in W \text{and any } (w R) x_0,\ldots,x_{n-1}$, there exists $(w R) z$, s.t. $z$ strongly unites $x_0,\ldots,x_{n-1}$ in $\mathfrak{F}$;
    \item $\mathfrak{F}$ satisfies $(su) \Leftrightarrow \text{for any } m \in \omega^*$, $\mathfrak{F}$ satisfies $(su_n)$;
\end{itemize}

\end{enumerate}
\end{definition}

\begin{lemma}\label{(su_1)}
For any  \textbf{S4} frame $\mathfrak{F} = \langle W, R \rangle$,
$\mathfrak{F}$ satisfies $(su_1)$.
\end{lemma}

\begin{proof}
For any  $w \in W$, any $(w R) x_0$, by the reflexivity of $R$, it is easy to prove that $x_0$ strongly unites $x_0$.
\qed  \end{proof}

\begin{lemma}[Strong union of strong unions]
For any  $n \in \omega^*$, any \textbf{S4} frame $\mathfrak{F} = \langle W, R \rangle$, any $x_0,\ldots,x_{n-1} \in W$, any $a, u \in W$,

\indent $u$ strongly unites $z, a$ and $z$ strongly unites $x_0,\ldots,x_{n-1} \Rightarrow u$ strongly unites $x_0,\ldots,x_{n-1}, a$.
\end{lemma}

\begin{proof}
On one hand, since $u$ strongly unites $z, a$, we have
\begin{align}
u R z, a \text{ and } u \in \Box(\blackdiamond\{z\} \cup \Diamond \blackdiamond\{a\}) \text{ and } u \in \Box(\blackdiamond\{a\} \cup \Diamond \blackdiamond\{z\}) \tag{*} \label{*}
\end{align}

On the other hand, since $z$ strongly unites $x_0,\ldots,x_{n-1}$, we have
\begin{align}
\text{for any } i \in n\text{, } z R x_i \text{ and } z \in \Box(\blackdiamond\{x_i\} \cup \bigcup_{i' \in n \setminus \{i\}}\Diamond \blackdiamond\{x_{i'}\})\text{.} \tag{**} \label{**}
\end{align}

Thus, for any  $i \in n$,
by $u R z, a$ (\eqref{*}), $z R x_0,\ldots,x_{n-1}$ (\eqref{*}), and the transitivity of $R$, we have $u R x_0,\ldots,x_{n-1}, a$, and:\\
\indent$u \in \Box(\blackdiamond\{z\} \cup \Diamond \blackdiamond\{a\})$\, $\cap$\, $\Box(\blackdiamond\{a\} \cup \Diamond \blackdiamond\{z\})$ \hfill\quad (\eqref{*})

\noindent$\Rightarrow u \in \Box(\blackdiamond\{x_{i}\} \cup \bigcup_{i' \in n \setminus \{i\}}\Diamond \blackdiamond\{x_{i'}\} \cup \Diamond \blackdiamond\{a\})$\, $\cap$\, $\Box(\blackdiamond\{a\} \cup \Diamond \blackdiamond\{x_{i}\} \cup$\\ $\bigcup_{i' \in n \setminus \{i\}}\Diamond \Diamond \blackdiamond\{x_{i'}\}$ $)$ \hfill\quad (\eqref{**}, $\Box$ monotonic, $\Diamond$ distributive over $\cup$)

\noindent$\Rightarrow u \in \Box(\blackdiamond\{x_{i}\} \cup \bigcup_{i' \in n \setminus \{i\}}\Diamond \blackdiamond\{x_{i'}\} \cup \Diamond \blackdiamond\{a\})$\, $\cap$\, $\Box(\blackdiamond\{a\} \cup \bigcup_{i \in n}\Diamond \blackdiamond\{x_{i} \})~~~~~~~~~~$ \hfil\quad ($R$ transitive)

\noindent$\Leftrightarrow u$ strongly unites $x_0,\ldots,x_{n-1}, a$ \hfill\quad (definition of strong union).
\qed  \end{proof}

\begin{lemma}[$(su_2) \Rightarrow (su_n)$]\label{(su_2) imply (su_n)}
For any \textbf{S4} frame $\mathfrak{F} = \langle W, R \rangle$,
\[\mathfrak{F} \text{ satisfies } (su_2) \ \Rightarrow\ \text{for any  } 2 \leq n \in \omega\text{, } \mathfrak{F} \text{ satisfies } (su_n)\text{.}\]
\end{lemma}

\begin{proof}
Induction on $n$.
The base case is obvious.

Let  $2 \leq n \in \omega$. Assume $\mathfrak{F}$ satisfies $(su_n)$.
Let $w \in W$, and $w R x_0,\ldots,x_n$.
On one hand, since $\mathfrak{F}$ satisfies $(su_n)$, there exists  $(w R) z$, $z$ strongly unites $x_0,\ldots,x_{n-1}$.
On the other hand, since $\mathfrak{F}$ satisfies $(su_2)$, there exists  $(w R) u$, $u$ strongly unites $z, x_n$.
Thus, by the previous lemma, $u$ strongly unites $x_0,\ldots,x_n$.
\qed  
\end{proof}

Combining Lemma \ref{(su_1)}, Lemma \ref{(su_2) imply (su_n)}, and the definition of (su), we have:

\begin{corollary}\label{(su) equal (su_2)}
For any  \textbf{S4} frame $\mathfrak{F} = \langle W, R \rangle$,
\[\mathfrak{F} \text{ satisfies } (su_2) \ \Leftrightarrow\ \mathfrak{F} \text{ satisfies } (su)\text{.}\]
\end{corollary}

Let $\,\vDash$ denote the satisfaction relation defined in the Kripke semantics of intuitionistic propositional logic \textbf{IPL} (please refer to, for example, \cite[p.~42]{gabbay_interpolation_2005} for details).

\begin{definition}\label{vDash_F,x, vDash_F}
    Let $\Gamma\subseteq Form$ and $\alpha\in Form$. For any \textbf{S4} frame $\mathfrak{F}=\left\langle W,R\right\rangle$ with $x\in W$, 
    define:
    \begin{align*}
        \Gamma \vDash_{\mathfrak{F},x} \alpha &\iff \textit{for any monotonic valuation } V ,\ ( \mathfrak{F},V,x \vDash \Gamma \Rightarrow \mathfrak{F},V,x\vDash \alpha);\\
        \Gamma \vDash_{\mathfrak{F}} \alpha &\iff \textit{for any }y\in W,\ \Gamma \vDash_{\mathfrak{F},y}\alpha.
    \end{align*} 
    ( $\mathfrak{F},V,x \vDash \Gamma$ means that $\mathfrak{F},V,x \vDash \phi$ for each $\phi\in\Gamma$.)\\
    If $\,\emptyset\vDash_{\mathfrak{F}} \alpha$, we simply write $\,\vDash_{\mathfrak{F}} \alpha$ and say that $\mathfrak{F}$ validates $\alpha$. 
\end{definition}

\begin{theorem}[Frame correspondence of $\boldsymbol{su}$]
For any  \textbf{S4} frame $\mathfrak{F} = \langle W, R \rangle$,
\[\mathfrak{F}\text{ validates } \boldsymbol{su} \ \Leftrightarrow\ \mathfrak{F} \text{ satisfies } (su_2)\text{.}\]
\end{theorem}

\begin{proof}
~
\begin{enumerate}

\item[``$\Leftarrow$'':] Assume $\mathfrak{F}$ satisfies $(su_2)$.

We prove that for any  monotonic valuation $V$ on $\mathfrak{F}$, $\Vert \boldsymbol{su} \Vert = W$, where $\Vert \boldsymbol{su} \Vert$ is the truth set of $\boldsymbol{su}$ in the model $\mathfrak{F},V$.

Let $A = \Vert \alpha \Vert, B = \Vert \beta \Vert, C = \Vert \varphi \Vert, D = \Vert \psi \Vert$, then:

\indent\quad~$\Vert ((\neg \varphi \rightarrow \psi) \land (\neg \psi \rightarrow \varphi)) \rightarrow \alpha \lor \beta \Rightarrow (\varphi \rightarrow \alpha) \lor (\psi \rightarrow \beta) \Vert = W$

$\Leftarrow \Vert ((\neg \varphi \rightarrow \psi) \land (\neg \psi \rightarrow \varphi)) \rightarrow \alpha \lor \beta \Vert \subseteq \Vert (\varphi \rightarrow \alpha) \lor (\psi \rightarrow \beta) \Vert$ \hfil\quad ($\Box W = W$, for any  $X, Y \subseteq W, ((W \setminus X) \cup Y = W \Leftrightarrow X \subseteq Y$))

$\Leftrightarrow \Diamond(C \cap -A) \cap \Diamond(D \cap -B) \subseteq \Diamond((\neg C \rightarrow D) \land (\neg D \rightarrow C) \cap -A \cap -B)$ \hfil\quad (definition of truth sets)

Now we prove that the inclusion in the last line holds. for any  $w \in W$, if $w$ belongs to the left-hand set, then there exist $(w R)x, y$, $(x \in C \cap -A$ and $y \in D \cap -B)$, and by $(su_2)$, there exists  $(w R)z$ such that $z$ strongly unites $x, y$ in $\mathfrak{F}$. Combining the fact that $A, B, C, D$ are $R$-upsets, we can prove that $z \in (\neg C \rightarrow D) \land (\neg D \rightarrow C) \cap -A \cap -B$.

\item[``$\Rightarrow$'':] We prove the contrapositive. Assume $\mathfrak{F}$ does not satisfy $(su_2)$. We prove that $\boldsymbol{su}$ is not valid on $\mathfrak{F}$.

Since $\mathfrak{F}$ does not satisfy $(su_2)$, there exists  $w \in W$, and $(w R)x, y$, such that for any  $(w R)z$,
\begin{align}
        z R x, y \Rightarrow z \in -\Box((\Diamond \blackdiamond\{x\} \cup \blackdiamond\{y\}) \cap (\Diamond \blackdiamond\{y\} \cup \blackdiamond\{x\}))  \tag{$\star$} \label{star}
\end{align}

Let $V = \{ \langle p, \blackdiamond\{x\} \rangle \} \cup \{\langle r, \neg\{x\} \rangle\} \cup \{\langle q, \blackdiamond\{y\} \rangle\} \cup \{\langle s, \neg\{y\} \rangle\} \cup \bigcup \{ \langle p, \emptyset \rangle | p \in PL \setminus \{r, s, p, q\} \}$.

On one hand, it is easy to prove that $\mathfrak{F}, V, x \vDash p$, $\mathfrak{F}, V, x \not\vDash r$, $\mathfrak{F}, V, y \vDash q$, $\mathfrak{F}, V, y \not\vDash s$, so $\mathfrak{F}, V, w \not\vDash (p \rightarrow r) \lor (q \rightarrow s)$.

On the other hand,

\indent\quad~$\mathfrak{F}, V, w \vDash ((\neg p \rightarrow q) \land (\neg q \rightarrow p)) \rightarrow r \lor s$

$\Leftrightarrow \text{for any } (w R)z, ( z \not\in \neg\{x\}$ and $z \not\in \neg\{y\} \Rightarrow z \in -\Box((- \neg \blackdiamond\{x\} \cup \blackdiamond\{y\}) \cap (- \neg \blackdiamond\{x\} \cup \blackdiamond\{y\})) )$ \hfill\quad (definition of $\vDash$, definition of $V$)

$\Leftrightarrow \text{for any } (w R)z, (z R x, y \Rightarrow z \in -\Box((\Diamond \blackdiamond\{x\} \cup \blackdiamond\{y\}) \cap (\Diamond \blackdiamond\{y\} \cup \blackdiamond\{x\})))$

$\Leftrightarrow$ truth \quad (\eqref{star}).

Thus, $\mathfrak{F}, V, w \not\vDash \boldsymbol{su}$ ($R$ is reflexive).
\end{enumerate}
\qed  
\end{proof}

\section{Syntactic and semantic consequence}
\begin{definition}[Syntactic consequence]
Let $\,\vdash_{\textbf{SU}}$ be the syntactic consequence relation obtained by adding $\boldsymbol{su}$ (as an axiom schema) to the Hilbert system of $\textbf{IPL}$.
\end{definition}

It is easy to see that $\,\vdash_{\textbf{IPL}}\ \subseteq\ \,\vdash_{\textbf{SU}}$. And as a superintuitionistic logic, $\vdash_{\textbf{SU}}$ has some basic properties, such as: (for simplicity, we omit the subscript of $\,\vdash$ and the universal quantifiers at the beginning of each clause)

\begin{enumerate}
  \item[(Com)] If $\Gamma\vdash\varphi$, then there is a finite $\Gamma'\subseteq\Gamma$ such that $\Gamma'\vdash\varphi$ ;
  \item[$(A)$]\ $\Gamma\cup\{\varphi\}\vdash\varphi$ ;
  \item[$(Cut)$] If $\Gamma\cup\{\psi\}\vdash\varphi$ and $\Delta\vdash\psi$, then $\Gamma\cup\Delta\vdash\varphi$ ;
  \item[(Mon)] If $\Gamma\vdash\varphi$ and $\Gamma\subseteq\Delta$, then $\Delta\vdash\varphi$ ;
  \item[$(\perp)$] $\perp \,\vdash\varphi$;
  \item[$(\land I)$] $\{\varphi,\psi\}\vdash\varphi\land\psi$ ;
  \item[$(\land E)$] $\varphi\land\psi\vdash\varphi$ and $\varphi\land\psi\vdash\psi$;
  \item[$(\vee I)$] $\varphi \vdash \varphi \vee \psi$ and $\psi \vdash \varphi \vee \psi$ ;
  \item[(PC)] If $\Gamma\cup\{\varphi\}\vdash\chi$ and $\Gamma\cup\{\psi\}\vdash\chi$, then $\Gamma\cup\{\varphi\vee\psi\}\vdash\chi$ ;
  \item[(MP)] $\{\varphi,\varphi\to\psi\}\vdash\psi$ ;
  \item[(DT)]  $\Gamma\cup\{\varphi\}\vdash\psi\Rightarrow\Gamma\vdash \varphi\to\psi$.
\end{enumerate}

Additionally, $\vdash_{\textbf{SU}}$ satisfies: 

\begin{lemma}
~
\begin{enumerate}
\item[(SU$^*$)]  $((\neg \varphi_{1} \rightarrow \psi_1) \land (\neg \psi_1 \rightarrow \varphi_{1}) \land \ldots \land (\neg \varphi_{n} \rightarrow \psi_n) \land (\neg \psi_n \rightarrow \varphi_{n})) \rightarrow \alpha \lor \beta \vdash_{\textbf{SU}} (\varphi_{1} \land \ldots \land \varphi_{n} \rightarrow \alpha) \lor (\psi_1 \land \ldots \land \psi_n \rightarrow \beta)$;

\end{enumerate}
\end{lemma}

\begin{proof}
By $\boldsymbol{su}$, we only need to prove $(\neg (\varphi_{1} \land \ldots \land \varphi_{n}) \rightarrow \psi_1 \land \ldots \land \psi_n) \land (\neg (\psi_1 \land \ldots \land \psi_n) \rightarrow \varphi_{1} \land \ldots \land \varphi_{n})$ $\vdash_{\textbf{SU}}$ $(\neg \varphi_{1} \rightarrow \psi_1) \land (\neg \psi_1 \rightarrow \varphi_{1}) \land \ldots \land (\neg \varphi_{n} \rightarrow \psi_n) \land (\neg \psi_n \rightarrow \varphi_{n})$. This is actually provable in intuitionistic logic.

\qed  
\end{proof}

\begin{definition}[Semantic consequence]
Let $\mathcal{C}_{su}$ be the class of all \textbf{S4} frames that satisfy $(su_2)$. We define the semantic consequence relation  of Strong Union logic \textbf{SU} as follows: for any  $\Gamma\subseteq Form$ and $\alpha\in Form$,
\begin{align*}
\Gamma \vDash_{\textbf{SU}} \alpha \ \Leftrightarrow\ \textit{for any }\mathfrak{F}\in \mathcal{C}_{su},\ \Gamma \vDash_{\mathfrak{F}}\alpha
\end{align*}
\end{definition}

Using the frame correspondence theorem, we can easily prove:

\begin{theorem}[Strong soundness theorem]\label{SU strong soundness}
For any  $\Gamma\subseteq Form$ and $\phi\in Form$,
\[\Gamma\,\vdash_{\textbf{SU}}\phi \ \Rightarrow\ \Gamma\,\vDash_{\textbf{SU}}\phi\textit{.}\]
\end{theorem}

\section{Completeness}
\subsection{Canonical model and its basic properties}
Let's recall some basic steps for proving Kripke completeness in intermediate logics: canonical model, Lindenbaum's lemma,  Truth lemma, etc.

\begin{definition}

Let $\ \,\vdash\ \subseteq \wp(Form) \times Form$. Define:

\begin{enumerate}

\item $\Gamma\in \wp(Form)$ is $\vdash$-closed if and only if for every $\varphi\in Form$, $\Gamma\vdash\varphi$ implies $\varphi\in\Gamma$.

\item $\Gamma\in \wp(Form)$ is $\vdash$-consistent if and only if $\Gamma\nvdash\perp$.

\item $\Gamma$ is $\,\vdash$-$\phi$-maximally consistent if and only if $\Gamma \,\nvdash \phi$ and for any  $\alpha \in Form \setminus \Gamma$, $\Gamma \cup \{\alpha\} \,\vdash\phi$.

\item $\Gamma\in \wp(Form)$ is disjunction complete if and only if for any  $\varphi,\psi\in Form$, if $\varphi\vee\psi\in\Gamma$, then $\varphi\in\Gamma$ or $\psi\in\Gamma$.

\end{enumerate}

\end{definition}

\begin{definition}

Let $\ \,\vdash\ \subseteq \wp(Form) \times Form$. Define

\begin{itemize}

\item $W^c_{\vdash}=\{\Gamma\subseteq Form\,|\,\ \Gamma$ is $\vdash$-consistent, $\vdash$-closed, and disjunction complete $\}$;

\item $V^c_{\vdash} : PL\to\wp(W^c_{\vdash})$::  $V^c_{\vdash} (p)=\{\Gamma\in W^c_{\vdash}\,|\, p\in\Gamma\}$ for each $p\in PL$.

\end{itemize}

Then define the canonical model relative to $\vdash$: $$\mathfrak{M}^c_{\vdash}=\langle\mathfrak{F}^c_{\vdash},V^c_{\vdash}\rangle=\langle W^c_{\vdash},\subseteq,V^c_{\vdash}\rangle.$$

Let $End_{\mathfrak{F}^c_{\vdash}}=\{\Phi\mid \Phi$ is a maximal element in $\mathfrak{F}^c_{\vdash}\}$, and for each $\Gamma\subseteq Form$, define $end(\Gamma)=\{\Phi\in W^c_{\vdash}\mid \Gamma\subseteq \Phi\}$.

\end{definition}

\begin{lemma}\label{superint: max.con. imply con. closed prime}

Let $\ \,\vdash\ \subseteq \wp(Form) \times Form$ satisfy $(\bot)$, (Cut), and (PC).
For any  $\Gamma \subseteq Form$, any $\phi \in Form$, if $\Gamma$ is $\,\vdash$-$\phi$-maximally consistent, then $\Gamma\in W^c_{\vdash}$.

\end{lemma}

\begin{proof}

Assume $\Gamma$ is $\,\vdash$-$\phi$-maximally consistent. The consistency and closedness of $\Gamma$ are straightforward. 
We now prove that $\Gamma$ is disjunction complete.

Let $\alpha_1 \vee \alpha_2 \in \Gamma$. By contradiction, assume $\alpha_1, \alpha_2 \notin \Gamma$. Since $\Gamma$ is $\,\vdash$-$\phi$-maximally consistent, $\Gamma \cup\{\alpha_i\} \vdash \phi$. Then, by (PC), we have $\Gamma \cup\{\alpha_1 \vee \alpha_2\}\vdash \phi$. Thus, $ \Gamma\vdash\phi$, which is a contradiction.
\qed  \end{proof}

\begin{lemma}[Lindenbaum's lemma]\label{superint: Lindenbaum lemma}
Let $\ \,\vdash\ \subseteq \wp(Form_\neg) \times Form_\neg$ satisfy \((Com), (Mon) \).
If \(\Gamma \nvdash\phi\), then there exists \(\Phi \subseteq \text{Form}_\neg\), such that \(\Phi\) is \(\,\vdash\)-\(\phi\)-maximal consistent and \(\Gamma \subseteq \Phi\) and $\phi\notin\Phi$.
\end{lemma}

\begin{proof}
Consider $\Omega = \{ \Delta \subseteq \text{Form}_\neg \mid \Gamma \subseteq \Delta \text{ and } \Delta \nvdash\phi \}$. 
It is easy to show that $\Gamma$ is in $\Omega$ (using $\Gamma\nvdash\phi$), and every chain in $\Omega$ has an upper bound (using (Com), (Mon)). 

By Zorn's lemma, $\Omega$ has a maximal element, denoted as $\Phi$. It's clear that $\Phi$ is $\vdash$-$\phi$-maximal consistent. 
\qed
\end{proof}

\begin{lemma}[$\rightarrow$-Existence lemma]\label{superint: existence lemma}
Let $\ \,\vdash\ \subseteq \wp(Form) \times Form$ satisfy (Cut), (Mon), (Com), $(\bot)$, (PC), (DT).
For any  $\Gamma \subseteq Form$ and $\alpha, \beta \in Form$, if $\Gamma \nvdash \alpha \rightarrow \beta$, then there exists  $\Phi \in W^c_{\vdash}$ such that $\Gamma \subseteq \Phi$ and $\alpha \in \Phi$ and $\beta \notin \Phi$.
\end{lemma}

\begin{proof}
Since $\Gamma \nvdash \alpha \rightarrow \beta$, by (DT), we know that $\Gamma\cup\{\alpha\} \nvdash \beta$. Then, using Lemma \ref{superint: Lindenbaum lemma} and Lemma \ref{superint: max.con. imply con. closed prime}, we can obtain the desired set $\Phi$.
\qed  \end{proof}

Combining the above lemmas, it is easy to prove by induction:

\begin{lemma}[Truth lemma]\label{superint: truth lemma}
Let $\ \,\vdash\ \subseteq \wp(Form) \times Form$ satisfy (A), (Cut), (Mon), (Com), $(\bot)$,  $(\land$I), $(\land$E), $(\lor$I), (PC), (DT), (MP).
For any  $\varphi\in Form$ and $\Gamma\in W^c_{\vdash}$,
\[\mathfrak{M}^c_{\vdash},\Gamma\vDash\varphi \Leftrightarrow \varphi\in\Gamma\textit{.}\]
\end{lemma}

\subsection{Canonical frame satisfies $(su_2)$}
Next, we give two crucial lemmas for demonstrating the completeness of $\mathbf{SU}$. Note that Lemma \ref{superint: canonical box(dia.RX cup RY)}, like those in the previous section, does not make use of $\boldsymbol{su}$.

\begin{lemma}\label{superint: canonical box(dia.RX cup RY)}
Let $\ \,\vdash\ \subseteq \wp(Form) \times Form$ satisfy (A), (Cut), (Mon), (Com), $(\bot)$, (PC), (MP), (DT).
For any  $\Gamma, \Delta, \Theta \subseteq \text{Form}$, if $\Gamma, \Delta, \Theta$ are all $\vdash$-closed, then
\begin{align*}
&\neg \Gamma \rightarrow \Delta = \{\neg \alpha \rightarrow \beta \ | \alpha \in \Gamma, \beta \in \Delta\} \subseteq \Theta\\
\Leftrightarrow &\text{for any } \Theta \subseteq \Theta' \in W^c_{\vdash}, (\Gamma \cup \Theta' \text{ is $\,\vdash$-consistent\ or\ } \Delta \subseteq \Theta')\\
\Leftrightarrow &\text{for any } \Theta \subseteq \Theta' \in W^c_{\vdash}, (\text{there exists } \Psi \in W, \Gamma \cup \Theta' \subseteq \Psi \text{\ or\ } \Delta \subseteq \Theta').
\end{align*}
\end{lemma}

\begin{proof}
~
\begin{itemize}

\item ``Top $\Rightarrow$ Middle'': Assume $\Gamma \cup \Theta'$ is not $\vdash$-consistent. By (Com), (Mon), and the fact that $\Gamma, \Delta$ are $\vdash$-closed, there exists  $\alpha \in \Gamma, \Theta' \cup \{\alpha\} \vdash \bot$. By (DT), $\Theta' \vdash \neg \alpha$.

Since $\neg \Gamma \rightarrow \Delta \subseteq \Theta \subseteq \Theta'$, by (A), for any  $\beta \in \Delta, \Theta' \vdash \neg \alpha \rightarrow \beta$. Thus, by (MP), (Cut), and the fact that $\Theta'$ is $\vdash$-closed, it is easy to see that $\Delta \subseteq \Theta'$.

\item ``Middle $\Rightarrow$ Top'': Assume there exists  $\alpha \in \Gamma,\beta \in \Delta, \neg \alpha \rightarrow \beta \notin \Theta$. Since $\Theta$ is $\vdash$-closed, $\Theta \nvdash \neg \alpha \rightarrow \beta$. By Lemma \ref{superint: existence lemma} (using $(\bot)$, (PC)), there exists  $\Phi \in W^c_{\vdash}, \Theta \cup \{\neg \alpha\} \subseteq \Phi$ and $\beta \notin \Phi$. This contradicts the "Middle" condition (using (MP)).

\item ``Middle $\Leftrightarrow$ Bottom'':
In fact, we can prove that for any  $\Gamma \subseteq \text{Form}$, ($\Gamma$ is $\vdash$-consistent $\Leftrightarrow$ there exists  $\Psi \in W^c_{\vdash}, \Gamma \subseteq \Psi$).
The "$\Rightarrow$" direction uses Lemma \ref{superint: Lindenbaum lemma} and Lemma \ref{superint: max.con. imply con. closed prime}, while the "$\Leftarrow$" direction follows easily from (Mon).
\end{itemize}
\qed  \end{proof}

\begin{lemma}\label{weak primeness of Psi_0}

Let $\Phi ,\Gamma, \Delta \in W^c_{\vdash_{\textbf{SU}}}$ such that $\Phi \subseteq \Gamma \cap \Delta$. Let $\Theta=\Phi \cup (\neg \Gamma \rightarrow \Delta) \cup (\neg \Delta \rightarrow \Gamma)$.  Then \( \Theta \subseteq \Gamma \cap \Delta\), and for any $n\in\mathbb{N}^*$ and any $\alpha_{0}, ..., \alpha_{n} \in \text{Form}$,
\begin{align}
\Theta \vdash_{\textbf{SU}} \alpha_{0} \lor...\lor \alpha_{n} \ \Rightarrow\ \alpha_{0} \in \Gamma \cap \Delta \text{\ or...or\ } \alpha_{n} \in \Gamma \cap \Delta \tag{$\#$} \label{prop.aa}
\end{align}

\end{lemma}

\begin{proof}

Using $\alpha \lor \beta \vdash_{\textbf{IPL}} \neg \alpha \rightarrow \beta$ and the fact that $\Gamma, \Delta$ are $\vdash_{\textbf{SU}}$-closed, it is easy to prove that $(\neg \Gamma \rightarrow \Delta) \cup (\neg \Delta \rightarrow \Gamma) \subseteq \Gamma \cap \Delta$. Combining with $\Phi \subseteq \Gamma \cap \Delta$, we have \( \Theta \subseteq \Gamma \cap \Delta\).

Next, we prove \eqref{prop.aa}.

Suppose $\Theta \vdash_{\textbf{SU}}\alpha_{0} \lor...\lor \alpha_{n}$. We want to prove that $\alpha_{0} \in \Gamma \cap \Delta \text{\ or\ ... or\ } \alpha_{n} \in \Gamma \cap \Delta$.

Using the distributive law, we only need to prove: $(\alpha_{0} \in \Gamma \text{\ or\ ...or\ } \alpha_{n} \in \Gamma) $ and $(\alpha_{0} \in \Delta \text{\ or\ ...or\ } \alpha_{n} \in \Delta)$ and for any permutation $\alpha_{i_0}, ..., \alpha_{i_n}$ of $\alpha_{0}, ..., \alpha_{n}$, and for any $0\leq m < n$, we have $(\alpha_{i_0} \in \Gamma \text{\ or\ ...or\ } \alpha_{i_m} \in \Gamma$ or $\alpha_{i_{m+1}} \in \Delta \text{\ or\ ...or\ } \alpha_{i_n} \in \Delta) $.

\begin{itemize}

\item From $\alpha_{0} \lor...\lor \alpha_{n} \in \Psi_{0}$, $\Phi \subseteq \Gamma \cap \Delta$, and the disjunction completeness of $\Gamma$ and $\Delta$, it is easy to see that $(\alpha_{0} \in \Gamma \text{\ or\ ...or\ } \alpha_{n} \in \Gamma) $ and $(\alpha_{0} \in \Delta \text{\ or\ ...or\ } \alpha_{n} \in \Delta)$.

\item For any permutation $\alpha_{i_{0}}, ..., \alpha_{i_{n}}$ of $\alpha_{0}, ..., \alpha_{n}$, and for any $0\leq m < n$, we have:

\indent\quad~$\Theta \vdash_{\textbf{SU}}\alpha_{0} \lor...\lor \alpha_{n}$

$\Rightarrow$ There exist $\varphi_{1}, \ldots, \varphi_{n} \in \Gamma$ and $\psi_1, \ldots, \psi_n \in \Delta$ such that $\Phi \vdash_{\textbf{SU}} ((\neg \varphi_{1} \rightarrow \psi_1 ) \land (\neg \psi_1 \rightarrow \varphi_{1}) \land \ldots \land (\neg \varphi_{n} \rightarrow \psi_n ) \land (\neg \psi_n \rightarrow \varphi_{n})) \rightarrow \alpha_{0} \lor...\lor \alpha_{n}$ \hfil\quad ((Com), (DT))

$\Leftrightarrow$ There exist $\varphi_{1}, \ldots, \varphi_{n} \in \Gamma$ and $\psi_1, \ldots, \psi_n \in \Delta$ such that $\Phi \vdash_{\textbf{SU}} ((\neg \varphi_{1} \rightarrow \psi_1) \land (\neg \psi_1 \rightarrow \varphi_{1}) \land \ldots \land (\neg \varphi_{n} \rightarrow \psi_n) \land (\neg \psi_n \rightarrow \varphi_{n})) \rightarrow (\alpha_{i_{0}} \lor...\lor \alpha_{i_{m}})\lor(\alpha_{i_{m+1}} \lor...\lor \alpha_{i_{n}})$ \hfil\quad ($\alpha_{i_0}, ..., \alpha_{i_n}$ is a permutation of $\alpha_{0}, ..., \alpha_{n}$)

$\Leftrightarrow$ There exist $\varphi_{1}, \ldots, \varphi_{n} \in \Gamma$ and $\psi_1, \ldots, \psi_n \in \Delta$ such that $\Phi \vdash_{\textbf{SU}} (\varphi_{1} \land \ldots \land \varphi_{n} \rightarrow \alpha_{i_{0}} \lor...\lor \alpha_{i_{m}}) \lor (\psi_1 \land \ldots \land \psi_n \rightarrow \alpha_{i_{m+1}} \lor...\lor \alpha_{i_{n}})$ \hfil\quad ((SU$^*$))

$\Leftrightarrow$ There exist $\varphi_{1}, \ldots, \varphi_{n} \in \Gamma$ and $\psi_1, \ldots, \psi_n \in \Delta$ such that $\Phi \vdash_{\textbf{SU}} \varphi_{1} \land \ldots \land \varphi_{n} \rightarrow \alpha_{i_{0}} \lor...\lor \alpha_{i_{m}} \text{\ or\ } \Phi \vdash_{\textbf{SU}} \psi_1 \land \ldots \land \psi_n \rightarrow \alpha_{i_{m+1}} \lor...\lor \alpha_{i_{n}}$ \hfil\quad ($\Phi$ is disjunction complete and $\vdash_{\textbf{SU}}$-closed)

$\Rightarrow \alpha_{i_{0}} \lor...\lor \alpha_{i_{m}} \in \Gamma \text{\ or\ } \alpha_{i_{m+1}} \lor...\lor \alpha_{i_{n}} \in \Delta$ \hfil\quad ($\Phi \subseteq \Gamma\cap \Delta$, (Mon), (MP), $\Gamma$ and $\Delta$ are $\vdash_{\textbf{SU}}$-closed)

$\Rightarrow$ $\alpha_{i_0} \in \Gamma \text{\ or\ ...or\ } \alpha_{i_m} \in \Gamma$ or $\alpha_{i_{m+1}} \in \Delta \text{\ or\ ...or\ } \alpha_{i_n} \in \Delta$ \hfil\quad ($\Gamma$ and $\Delta$ are disjunction complete).
\end{itemize}
\qed
\end{proof}

Now, we are at the point to prove that the $\,\vdash_{SU}$-canonical frame satisfies ($su_{2}$). 

{\ttfamily For any $\Phi ,\Gamma, \Delta \in W^c_{\vdash_{\textbf{SU}}}$, if $\Phi \subseteq \Gamma \cap \Delta$, then there exists $\Psi \in W^c_{\vdash_{\textbf{SU}}}$,
such that $\Phi \subseteq \Psi \subseteq \Gamma \cap \Delta$ and $\Psi \in \Box^c((\Diamond^c \blackdiamond^c\{\Gamma\} \cup \blackdiamond^c\{\Delta\}) \cap (\Diamond^c \blackdiamond^c\{\Delta\} \cup \blackdiamond^c\{\Gamma\}))$.}

(Here, $\blackdiamond^c$, $\Box^c$, and $\Diamond^c$ stand for the operators $\blackdiamond^c_{\vdash_{\textbf{SU}}}$, $\Box_{\mathfrak{F^c_{\vdash_{\textbf{SU}}}}}$, and $\Diamond_{\mathfrak{F^c_{\vdash_{\textbf{SU}}}}}$ on $\wp(W^c_{\vdash_{\textbf{SU}}})$.)

By Lemma \ref{superint: canonical box(dia.RX cup RY)}, we only need to construct, for any $\Phi ,\Gamma, \Delta \in W^c_{\vdash_{\textbf{SU}}}$ with $\Phi \subseteq \Gamma \cap \Delta$, a $\Psi \in W^c_{\vdash_{\textbf{SU}}}$ such that $\Phi\cup (\neg \Gamma \rightarrow \Delta) \cup (\neg \Delta \rightarrow \Gamma) \subseteq \Psi \subseteq \Gamma \cap \Delta$.

Our approach is to list all disjunctive formulas, then starting with $\Psi_{0}=\Phi \cup (\neg \Gamma \rightarrow \Delta) \cup (\neg \Delta \rightarrow \Gamma)$, construct a sequence $\Psi_{0}\subseteq\Psi_{1}\subseteq\Psi_{2}\subseteq$... At the $(k+1)$-th step, we handle the first disjunctive formula that has not been processed yet and is derivable from $\Psi_{k}$, adding one of its disjuncts to $\Psi_{k}$. Finally, we take $\Psi=\bigcup\Psi_{i}$.

However, the key is to ensure that at each step, we can always find a disjunct that belongs to $\Gamma \cap \Delta$, thereby guaranteeing that each $\Psi_{k}$ is contained within $\Gamma \cap \Delta$.

Let us represent all possible construction paths using a binary tree. (In this way, an ``good" infinite construction path corresponds to an infinite branch of a specific subtree (the subtree consisting of all ``good" points) of the binary tree. The existence of such an infinite branch, by König's lemma, can be reduced to the infiniteness of the subtree.)

\begin{definition}

Let $\emptyset$ be the unique 01-sequence of length 0. For any finite 01-sequences $\bar{s},\bar{t}$, define:

\begin{enumerate}

\item $\bar{s}\bar{t}$ is the 01-sequence obtained by concatenating $\bar{t}$ after $\bar{s}$.

\item $\bar{t}\preceq\bar{s}$ means that $\bar{t}$ is a initial segment of $\bar{s}$.

\end{enumerate}

\end{definition}

\begin{figure}[!ht]
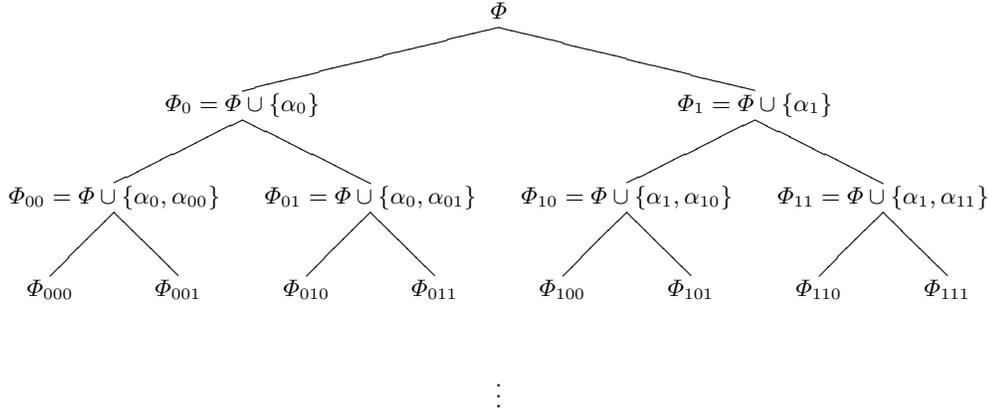


\centering

\xytree[0.02]{
            &  &  & & & & & \xynode[-4,4]{$\Phi$} & & & &  &  & &\\
            & &  & \xynode[-2,2]{$\Phi_{0}=\Phi\cup\{\alpha_{0}\}$} &  & &  &  & &  & & \xynode[-2,2]{$\Phi_{1}=\Phi\cup\{\alpha_{1}\}$} &  & & \\
           & \xynode[-1,1]{$\Phi_{00}=\Phi\cup\{\alpha_{0},\alpha_{00}\}$} & & & & \xynode[-1,1]{$\Phi_{01}=\Phi\cup\{\alpha_{0},\alpha_{01}\}$} & & &
           & \xynode[-1,1]{$\Phi_{10}=\Phi\cup\{\alpha_{1},\alpha_{10}\}$} & & & & \xynode[-1,1]{$\Phi_{11}=\Phi\cup\{\alpha_{1},\alpha_{11}\}$} &\\
           \xynode{$\Phi_{000}$} &  & \xynode{$\Phi_{001}$} & & \xynode{$\Phi_{010}$} & & \xynode{$\Phi_{011}$} & &
           \xynode{$\Phi_{100}$} & & \xynode{$\Phi_{101}$} & & \xynode{$\Phi_{110}$} &  & \xynode{$\Phi_{111}$}\\
           &  &  & & & & & \xynode{$\vdots$} & & & &  &  & &
            }

\caption{$\Phi$-disjunction tree}

\end{figure}

\begin{definition}[$\Phi$-disjunction tree]
Let all disjunctive formulas be enumerated as $(\phi_{0}\lor\psi_{0}), (\phi_{1}\lor\psi_{1}), $... Denote this sequence of formulas as $\Xi$. For any set of formulas $\Phi$, recursively define the following labeling on the infinite binary tree (i.e., the partial order structure of all finite 01-sequences under the $\preceq$ relation):

\begin{itemize}

\item $\alpha_{\emptyset}=T$, $\Phi_{\emptyset}=\Phi$, $\Xi_{\emptyset}=\Xi$;

\item For any finite 01-sequence $\bar{s}$, let $\alpha_{\bar{s}0}\lor \alpha_{\bar{s}1}$ be the first disjunctive formula in $\Xi_{\bar{s}}$ that can be derived from $\Phi_{\bar{s}}$, and then let $\Phi_{\bar{s}0} = \Phi_{\bar{s}}\cup\{\alpha_{\bar{s}0}\}$, $\Phi_{\bar{s}1} = \Phi_{\bar{s}}\cup\{\alpha_{\bar{s}1}\}$, and $\Xi_{\bar{s}0} = \Xi_{\bar{s}1}$ be the sequence of disjunctive formulas obtained by removing $\alpha_{\bar{s}0}\lor \alpha_{\bar{s}1}$ from $\Xi_{\bar{s}}$.

\end{itemize}

\end{definition}

Here is a basic property of the disjunction tree:

\begin{lemma}\label{Disjuntion Tree: full disjuntion span is provable}

For any $k\in\mathbb{N}$, let $S_{k}$ be the set of all 01-sequences of length $k$. For any set of formulas $\Phi$, and any point $\bar{r}$ in the $\Phi$-disjunction tree, if $f$ is a function on $S_{k}$ such that for any $\bar{s}\in S_{k}$, $f(\bar{s})\preceq \bar{s}$, then

\[\Phi_{\bar{r}}\vdash_{\textbf{SU}}\bigvee\{\alpha_{\bar{r}f(\bar{s})}\mid \bar{s}\in S_{k}\}.\]

\end{lemma}

\begin{proof}

By induction on $k$.

When $k=0$, $S_{0}=\{\emptyset\}$, and it is easy to see that $f(\emptyset)=\emptyset$. By the definition of the disjunction tree, $\alpha_{\bar{r}}\in\Phi_{\bar{r}}$, so $\Phi_{\bar{r}}\vdash_{\textbf{SU}}\bigvee\{\alpha_{\bar{r}}\}$.

Assume the proposition holds for $k$. Let $f$ be a function on $S_{k+1}$ such that for any $\bar{u}\in S_{k+1}$, $f(\bar{u})\preceq \bar{u}$.

The sequences in $S_{k+1}$ can be divided into two categories based on whether the first bit is 0 or 1, resulting in a partition of $S_{k+1}$ into $0S_{k}$ and $1S_{k}$:

\begin{itemize}

\item[] $0S_{k}=\{0\bar{s}\mid \bar{s}\in S_{k}\}$,

\item[] $1S_{k}=\{1\bar{s}\mid \bar{s}\in S_{k}\}$.

\end{itemize}

We only need to prove that $\Phi_{\bar{r}0}\vdash_{\textbf{SU}}\bigvee\{\alpha_{\bar{r}f(\bar{u})}\mid \bar{u}\in 0S_{k}\}$ and $\Phi_{\bar{r}1}\vdash_{\textbf{SU}}\bigvee\{\alpha_{\bar{r}f(\bar{u})}\mid \bar{u}\in 1S_{k}\}$. (------ If proven, then by ($\lor$I), both $\Phi_{\bar{r}0}$ and $\Phi_{\bar{r}1}$ can derive $\bigvee\{\alpha_{\bar{r}f(\bar{u})}\mid \bar{u}\in S_{k+1}\}$, and combined with (PC), we have $\Phi_{\bar{r}}\cup\{\alpha_{\bar{r}0}\lor\alpha_{\bar{r}1}\}\vdash_{\textbf{SU}}\bigvee\{\alpha_{\bar{r}f(\bar{u})}\mid \bar{u}\in S_{k+1}\}$. Since by the definition of the disjunction tree, $\Theta_{\bar{r}}\vdash_{\textbf{SU}}\alpha_{\bar{r}0}\lor\alpha_{\bar{r}1}$, by (Cut), we have $\Phi_{\bar{r}}\vdash_{\textbf{SU}}\bigvee\{\alpha_{\bar{r}f(\bar{u})}\mid \bar{u}\in S_{k+1}\}$.)

\begin{itemize}

\item If there exists $\bar{u}\in 0S_{k}$ such that $f(\bar{u})=\emptyset$, then $\bigvee\{\alpha_{\bar{r}f(\bar{u})}\mid \bar{u}\in 0S_{k}\}$ has a disjunct $\alpha_{\bar{r}}$, and since $\alpha_{\bar{r}}\in\Phi_{\bar{r}}\subseteq\Phi_{\bar{r}0}$, it is easy to prove that $\Phi_{\bar{r}0}\vdash_{\textbf{SU}}\bigvee\{\alpha_{\bar{r}f(\bar{u})}\mid \bar{u}\in 0S_{k}\}$.

Otherwise, $f$ is non-empty on $0S_{k}$, and since $f$ outputs a initial segment of the input, there exists a function $g$ on $S_{k}$ such that for any $\bar{s}\in S_{k}$, $g(\bar{s})\preceq \bar{s}$, and for any $\bar{s}\in S_{k}$, $f(0\bar{s})= 0g(\bar{s})$. Applying the induction hypothesis to $\bar{r}0$ and $g$, we have $\Phi_{\bar{r}0}\vdash_{\textbf{SU}}\bigvee\{\alpha_{\bar{r}0g(\bar{s})}\mid \bar{s}\in S_{k}\}$, and thus $\Phi_{\bar{r}0}\vdash_{\textbf{SU}}\bigvee\{\alpha_{\bar{r}f(\bar{u})}\mid \bar{u}\in 0S_{k}\}$.

\item The proof for $\Phi_{\bar{r}1}\vdash_{\textbf{SU}}\bigvee\{\alpha_{\bar{r}f(\bar{u})}\mid \bar{u}\in 1S_{k}\}$ is similar.

\end{itemize}

\end{proof}

\begin{lemma}[Infiniteness of the desired subtree]\label{infiniteness of subtree}

Let $\Phi ,\Gamma, \Delta \in W^c_{\vdash_{\textbf{SU}}}$ satisfy $\Phi \subseteq \Gamma \cap \Delta$. Let $\Theta=\Phi \cup (\neg \Gamma \rightarrow \Delta) \cup (\neg \Delta \rightarrow \Gamma)$. Then, for any $k\in\mathbb{N}$, there exists a point $\bar{s}$ in the $k$-th layer of the $\Theta$-disjunction tree such that $\Theta_{\bar{s}}\subseteq\Gamma \cap \Delta$.

\end{lemma}

\begin{proof}
Let $\Phi ,\Gamma, \Delta \in W^c_{\vdash_{\textbf{SU}}}$ satisfy $\Phi \subseteq \Gamma \cap \Delta$. Let $\Theta=\Phi \cup (\neg \Gamma \rightarrow \Delta) \cup (\neg \Delta \rightarrow \Gamma)$. For any $k\in\mathbb{N}$, let the points in the $k$th layer of the $\Theta$-disjunction tree be $\bar{s}_{1},...,\bar{s}_{m}$. We want to prove that there exists $1\leq i\leq m$ such that $\Theta_{\bar{s}_{i}}\subseteq\Gamma \cap \Delta$.

By the definition of $\Theta_{\bar{s}_{i}}$, the goal is equivalent to the existence of $1\leq i\leq m$ such that $\{\alpha_{\bar{t}} \mid \bar{t}\preceq\bar{s}_{i}\}\subseteq\Gamma \cap \Delta$. This is actually a disjunction of $m$ conjunctive propositions. Using the distributive law, this is equivalent to a conjunction of some disjunctions. Specifically, after applying the distributive law, we only need to prove: for any $\bar{t}_{1},...,\bar{t}_{m}$, if $\bar{t}_{i}$ are initial segments of $\bar{s}_{i}$ respectively, then $\alpha_{\bar{t}_{1}}\in\Gamma \cap \Delta$ or...or $\alpha_{\bar{t}_{m}}\in\Gamma \cap \Delta$.

Since $\bar{t}_{i}$ are initial segments of $\bar{s}_{i}$ respectively, by Lemma \ref{Disjuntion Tree: full disjuntion span is provable}, we have $\Theta\vdash_{\textbf{SU}}\alpha_{\bar{t}_{1}}\lor...\lor\alpha_{\bar{t}_{m}}$. Thus, combining with Lemma \ref{weak primeness of Psi_0}, we have $\alpha_{\bar{t}_{1}}\in\Gamma \cap \Delta$ or...or $\alpha_{\bar{t}_{m}}\in\Gamma \cap \Delta$.

\end{proof}

From Lemma \ref{infiniteness of subtree} and König's lemma, we have:

\begin{corollary}

If $\Phi ,\Gamma, \Delta \in W^c_{\vdash_{\textbf{SU}}}$ satisfy $\Phi \subseteq \Gamma \cap \Delta$, let $\Theta=\Phi \cup (\neg \Gamma \rightarrow \Delta) \cup (\neg \Delta \rightarrow \Gamma)$, then the subtree of the $\Theta$-disjunction tree on $\{\bar{s}\mid\Theta_{\bar{s}}\subseteq\Gamma \cap \Delta\}$ has an infinite branch.

\end{corollary}

Now we can prove that the canonical frame $\mathfrak{F}^c_{\vdash_{\textbf{SU}}}$ satisfies (su$_2$):

\begin{lemma}\label{canonical frame satisfy (su_2)}

For any $\Phi ,\Gamma, \Delta \in W^c_{\vdash_{\textbf{SU}}}$, if $\Phi \subseteq \Gamma \cap \Delta$, then there exists $\Psi \in W^c_{\vdash_{\textbf{SU}}}$,

such that $\Phi \subseteq \Psi \subseteq \Gamma \cap \Delta$ and $\Psi \in \Box^c((\Diamond^c \blackdiamond^c\{\Gamma\} \cup \blackdiamond^c\{\Delta\}) \cap (\Diamond^c \blackdiamond^c\{\Delta\} \cup \blackdiamond^c\{\Gamma\}))$.

\end{lemma}

\begin{proof}

Let $\Theta=\Phi \cup (\neg \Gamma \rightarrow \Delta) \cup (\neg \Delta \rightarrow \Gamma)$. By Lemma \ref{superint: canonical box(dia.RX cup RY)}, we only need to find a $\Psi \in W^c_{\vdash_{\textbf{SU}}}$ such that $\Theta \subseteq \Psi \subseteq \Gamma \cap \Delta$.

From the above corollary, let $(\emptyset=)\bar{s}_{0},\bar{s}_{1},\bar{s}_{2}...$ be an infinite branch of the subtree of the $\Theta$-disjunction tree on $\{\bar{s}\mid\Theta_{\bar{s}}\subseteq\Gamma \cap \Delta\}$. Let $\Psi=\bigcup\Theta_{\bar{s}_{i}}$. It is easy to see that $\Theta\subseteq\Psi\subseteq\Gamma \cap \Delta$.

Below, we prove that $\Psi\in W^c_{\vdash_{\textbf{SU}}}$.

\begin{itemize}

\item From $\Psi\subseteq\Gamma \cap \Delta$ and the $\,\vdash_{\textbf{SU}}$-consistency of $\Gamma$, it is easy to prove the $\,\vdash_{\textbf{SU}}$-consistency of $\Psi$.

\item We prove a slightly stronger proposition than "disjunction complete": for any $\phi,\psi\in Form$, if $\Psi\vdash_{\textbf{SU}}\phi\lor\psi$, then $\phi\in\Psi$ or $\psi\in\Psi$.

Suppose $\Psi\vdash_{\textbf{SU}}\phi\lor\psi$. Then by (Com), there exists $i\in\mathbb{N}$ such that $\Theta_{\bar{s}_{i}}\vdash_{\textbf{SU}}\phi\lor\psi$. Since we enumerated all disjunctive formulas when defining the $\Theta$-disjunction tree, $\phi\lor\psi$ will appear in the enumeration, and there are only finitely many disjunctive formulas before it, say $j$. Since the disjunction tree handles the ``selection problem'' of a new disjunctive formula at each branching, by the time we reach $\Theta_{\bar{s}_{i+j+1}}$, the selection of $\phi\lor\psi$ has already been processed, i.e., either $\phi\in\Theta_{\bar{s}_{i+j+1}}$ or $\psi\in\Theta_{\bar{s}_{i+j+1}}$. Thus, $\phi\in\Psi$ or $\psi\in\Psi$.

\item For any $\phi\in Form$, if $\Psi\vdash_{\textbf{SU}}\phi$, then $\Psi\vdash_{\textbf{SU}}\phi\lor\phi$, and by the previous result, we have $\phi\in\Psi$.

\end{itemize}

\end{proof}

Thus, combining Lemma \ref{canonical frame satisfy (su_2)} with Lemma \ref{superint: max.con. imply con. closed prime}, Lemma \ref{superint: Lindenbaum lemma}, , and Lemma \ref{superint: truth lemma}, it is easy to prove,

\begin{theorem}[Strong Completeness Theorem]\label{SU strong completeness}
For any $\Gamma\subseteq Form$ and $\phi\in Form$,

\[\Gamma\,\vDash_{\textbf{SU}}\phi \ \Rightarrow\ \Gamma\,\vdash_{\textbf{SU}}\phi\textit{.}\]
\end{theorem}

\section{Disjunction Property}

\begin{definition}[Hereditary union function]
For any  \textbf{S4} frame $\mathfrak{F} = \langle W, R \rangle$, and any $f: W^2(\text{par}) \rightarrow W$, we say that $f$ is a ``hereditary union function'' in $\mathfrak{F}$, if and only if, $\text{dom}(f)$ is closed under taking $R^\sharp$-successors, and for any  $\langle w, v \rangle \in \text{dom}(f)$, we have: 

$f(w, v) R w, v$, and for any  $(f(w, v)R)t$, $(wRt \text{ or } vRt$ or there exists  $(wR)w'$ and $(vR)v'$ with $t=f(w', v'))$.\\
(Here, $R^\sharp = \{\langle \langle w, v \rangle, \langle w', v' \rangle \rangle \in W^2 \times W^2 \mid wRw' \text{\ and\ } vRv'\}$.)

Moreover, $f$ is called normal, if and only if, for any  $\langle w, v \rangle, \langle w', v' \rangle \in \text{dom}(f)$, $(wRw'$ and $vRv'\ \Rightarrow\ f(w, v)Rf(w', v'))$.
\end{definition}

\begin{lemma}[Hereditary union function outputs strong union]\label{-2. Hereditary Union Function Yields Strong Union}
For any  $\langle w, v \rangle \in \text{dom}(f)$,
$f(w, v)$ strongly unites $w, v$ in $\mathfrak{F}$.
\end{lemma}

\begin{proof}
~
\begin{itemize}
    \item $f(w, v) R w, v$: This follows directly from the definition of the hereditary union function.
    \item $f(w, v) \in \Box(\Diamond \blackdiamond\{w\} \cup \blackdiamond\{v\}) \cap \Box(\Diamond \blackdiamond\{v\} \cup \blackdiamond\{w\})$:
    We only prove $f(w, v) \in \Box(\Diamond \blackdiamond\{w\} \cup \blackdiamond\{v\})$, the other part is similar.
    For any  $(f(w, v)R)t$, by the definition of the hereditary union function, we consider the following cases:
    If $wRt$, then by the reflexivity of $R$, it is easy to prove that $t \in \Diamond \blackdiamond\{w\}$;
    If $vRt$, then $t \in \blackdiamond\{v\}$;
    Otherwise, there exists  $(wR)w'$, $(vR)v', t=f(w', v')$, then by $f(w', v')Rw'$, $wRw'$, we have $t \in \Diamond \blackdiamond\{w\}$.
\end{itemize}
\qed  \end{proof}

\begin{lemma}\label{-1. Hereditary Union Function and Strong Union}
For any  \textbf{S4} frame $\mathfrak{F} = \langle W, R \rangle$, and any normal hereditary union function $f$ in $\mathfrak{F}$,
for any  $m, n, l \in \omega$, any $a_0, \ldots, a_{m-1}, b_0, \ldots, b_{n-1} \in W$, any $\langle x_0, y_0 \rangle, \ldots, \langle x_{l-1}, y_{l-1} \rangle \in \text{dom}(f)$, and any $u_1, u_2 \in W$, if $m+l \geq 1$ and $n+l \geq 1$, then

$u_1$ strongly unites $a_0, \ldots, a_{m-1}, x_0, \ldots, x_{l-1}$ while $u_2$ strongly unites $b_0, \ldots, b_{n-1},$ $y_0, \ldots, y_{l-1}$ \\
$\Longrightarrow$ $f(u_1, u_2)$ strongly unites $a_0, \ldots, a_{m-1}, f(x_0, y_0), \ldots, f(x_{l-1}, y_{l-1}), b_0, \ldots, b_{n-1}$.
\end{lemma}

\begin{proof}
~
\begin{itemize}

\item $f(u_1, u_2) R f(x_i, y_i)$: This follows from the fact that $f$ is normal, $u_1 R x_i$, and $u_2 R y_i$.

\item $f(u_1, u_2) R a_j, b_k$: This follows from $f(u_1, u_2) R u_1, u_2$, $u_1 R a_j$, $u_2 R b_k$, and the transitivity of $R$.

\item We prove: for any  $i \in l$, $f(u_1, u_2) \in \Box(\blackdiamond\{f(x_i, y_i)\} \cup \bigcup_{i' \in l\setminus\{i\}}\Diamond \blackdiamond\{f(x_{i'}, y_{i'})\}$ $\cup$ $\bigcup_{j \in m}\Diamond \blackdiamond\{a_j\} \cup \bigcup_{k \in n}\Diamond \blackdiamond\{b_k\})$.

Let  $i \in l$ and $(f(u_1, u_2)R)w$.
We need to prove that $w \in \blackdiamond\{f(x_i, y_i)\} \cup \bigcup_{i' \in l\setminus\{i\}}\Diamond \blackdiamond\{f(x_{i'}, y_{i'})\} \cup \bigcup_{j \in m}\Diamond \blackdiamond\{a_j\} \cup \bigcup_{k \in n}\Diamond \blackdiamond\{b_k\}$.

By the definition of the hereditary union function, we consider the following cases:
\begin{itemize}
    \item If $u_1 R w$, then by the strong union property of $u_1$, we have $w \in \blackdiamond\{x_i\} \cup$ $\bigcup_{i' \in l\setminus\{i\}}\Diamond \blackdiamond\{x_{i'}\} \cup \bigcup_{j \in m}\Diamond \blackdiamond\{a_j\}$. Since $f(x_i, y_i)R x_i$ and $R$ is transitive, we have $\blackdiamond\{x_i\} \subseteq \blackdiamond\{f(x_i, y_i)\}$. Thus, by the monotonicity of the operators, the goal is easily proved.
    \item If $u_2 R w$, then a similar argument as above can be made.
    \item Otherwise, there exist $(u_1R)t_1$ and $(u_2R)t_2$ with $w=f(t_1, t_2)$. By the strong union properties of $u_1$ and $u_2$, we have $t_1 \in \blackdiamond\{x_i\} \cup \bigcup_{i' \in l\setminus\{i\}}\Diamond R$ $\{x_{i'}\} \cup \bigcup_{j \in m}\Diamond \blackdiamond\{a_j\}$ and $t_2 \in \blackdiamond\{y_i\} \cup \bigcup_{i' \in l\setminus\{i\}}\Diamond \blackdiamond\{y_{i'}\} \cup \bigcup_{k \in n}\Diamond \blackdiamond\{b_k\}$ -- denoted as (*) . Then

\indent\quad~$f(t_1, t_2) \in \blackdiamond\{f(x_i, y_i)\} \cup \bigcup_{i' \in l\setminus\{i\}}\Diamond \blackdiamond\{f(x_{i'}, y_{i'})\} \cup \bigcup_{j \in m}\Diamond \blackdiamond\{a_j\} \cup$\\ $\bigcup_{k \in n}\Diamond \blackdiamond\{b_k\}$

$\Leftarrow (t_1 \in \blackdiamond\{x_i\} \text{\ and\ } t_2 \in \blackdiamond\{y_i\}) \text{\ or\ } (t_1 \in \bigcup_{i' \in l\setminus\{i\}}\Diamond \blackdiamond\{x_{i'}\} \cup \bigcup_{j \in m}\Diamond \blackdiamond\{a_j\}$ or $t_2 \in \bigcup_{i' \in l\setminus\{i\}}\Diamond \blackdiamond\{y_{i'}\} \cup \bigcup_{k \in n}\Diamond \blackdiamond\{b_k\})$ \hfil\quad (The "$\Leftarrow$" in the left branch: $f$ is normal; the "$\Leftarrow$" in the right branch: $f(t_1, t_2)R t_1, t_2$, $f(x_{i'}, y_{i'})R x_{i'}, y_{i'}$, and $R$ is transitive.)

This holds (using (*)).
\end{itemize}

\item We now prove: for any  $j \in m$, $f(u_1, u_2) \in \Box(\blackdiamond\{a_j\} \cup \bigcup_{j' \in m\setminus\{j\}}\Diamond \blackdiamond\{a_{j'}\} \cup \bigcup_{i \in l}\Diamond \blackdiamond\{f(x_i, y_i)\} \cup \bigcup_{k \in n}\Diamond \blackdiamond\{b_k\})$.

Let  $j \in m$ and $(f(u_1, u_2)R)w$.

We need to prove that $w \in \blackdiamond\{a_j\} \cup \bigcup_{j' \in m\setminus\{j\}}\Diamond \blackdiamond\{a_{j'}\} \cup \bigcup_{i \in l}\Diamond \blackdiamond\{f(x_i, y_i)\} \cup \bigcup_{k \in n}\Diamond \blackdiamond\{b_k\}$.

By the definition of the hereditary union function, we consider the following cases:
\begin{itemize}
    \item If $u_1 R w$, then by the strong union property of $u_1$, we have $w \in \blackdiamond\{a_j\} \cup \bigcup_{j' \in m\setminus\{j\}}\Diamond \blackdiamond\{a_{j'}\} \cup \bigcup_{i \in l}\Diamond \blackdiamond\{x_i\}$. Since $f(x_i, y_i)R x_i$ and $R$ is transitive, we have $\blackdiamond\{x_i\} \subseteq \blackdiamond\{f(x_i, y_i)\}$. Thus, by the monotonicity of the operators, the goal is easily proved;
    \item If $u_2 R w$, then a similar argument as above can be made;
    \item Otherwise, there exist $(u_1R)t_1$ and $(u_2R)t_2$ with $w=f(t_1, t_2)$. By $u_2 R t_2$, the strong union property of $u_2$, and the reflexivity of $R$, we have $t_2 \in \bigcup_{k \in n}\Diamond \blackdiamond\{b_k\} \cup \bigcup_{i \in l}\Diamond \blackdiamond\{y_i\}$. 
    Meanwhile,  by $f(t_1, t_2)R t_2$, $f(x_i, y_i)R y_{i'}$, and the transitivity of $R$, we have $f(t_1, t_2) \in \bigcup_{k \in n}\Diamond \blackdiamond\{b_k\} \cup$ $\bigcup_{i \in l}\Diamond \blackdiamond$ $\{f(x_i, y_i)\}$, so the goal is easily achieved.
\end{itemize}

\item Similarly one can prove: for any  $k \in n$, $f(u_1, u_2) \in \Box(\blackdiamond\{b_k\} \cup \bigcup_{k' \in n\setminus\{k\}}\Diamond \blackdiamond\{b_{k'}\}$ $\cup$ $\bigcup_{i \in l}\Diamond \blackdiamond\{f(x_i, y_i)\} \cup \bigcup_{j \in m}\Diamond \blackdiamond\{a_j\})$.

\end{itemize}
\qed  \end{proof}

Now we introduce a method to merge frames, which has been first considered by G. C. Meloni \cite{meloni_modelli_1984}.
\begin{definition}[connected product]
For any  \textbf{S4} frames $\mathfrak{F}_1 = \langle W_1, R_1 \rangle$ and $\mathfrak{F}_2 = \langle W_2, R_2 \rangle$, where $W_1, W_2$ are disjoint.
Define the connected product (following the terminology in \cite[p.~164]{vio_unification_2023}) of $\mathfrak{F}_1$ and $\mathfrak{F}_2$ as $\mathfrak{F}_1 \otimes \mathfrak{F}_2=\langle W, R \rangle$:

$W = W_1 \cup W_2 \cup (W_1 \times W_2)$,

$R = R_1 \cup R_2 \cup \{\langle \langle w_1, w_2\rangle , t \rangle \mid w_1R_1t \text{ or } w_2R_2t \text{ or } \text{there exists } (w_1R_1)v_1$ and $(w_2R_2)v_2$ with $t=\langle v_1, v_2\rangle \}$.
\end{definition}

\begin{lemma}[Basic properties of connected product]\label{0. HU-merging Frame Basic Properties}
Let $\mathfrak{F}=\mathfrak{F}_1 \otimes \mathfrak{F}_2=\langle W, R \rangle$. (For simplicity, we omit the universal quantifier ``for any  $i \in \{1,2\}$'' at the beginning of each sentence.)

\begin{enumerate}
\item 
  \begin{itemize}
    \item For any  $w_i \in W_i$, $\blackdiamond\{w_i\}=\blackdiamond_i\{w_i\}$.

    \item For any  $w_i, v_i \in W_i$, $(\langle w_1, w_2\rangle R v_i \Leftrightarrow w_i R_i v_i)$, and $(\langle w_1, w_2\rangle R$ $\langle v_1, v_2\rangle $ $\Leftrightarrow$ $w_1 R_1 v_1 \text{\ and\ } w_2 R_2 v_2)$.

    \item For any  $w \in W$, any $v_1 \in W_1$, any $v_2 \in W_2$, $(wR v_1, v_2 \Rightarrow \text{ there exists }$ $ \langle w_1, w_2 \rangle \in$ $W_1 \times W_2$, s.t. $(w_1 R_1 v_1 \text{\ and\ } w_2 R_2 v_2 \text{\ and\ } w=\langle w_1, w_2\rangle ))$.
  \end{itemize}

\item $R$ is reflexive and transitive.

\item If $r_1, r_2$ are the roots of $\mathfrak{F}_1, \mathfrak{F}_2$ respectively, then $\langle r_1, r_2\rangle $ is the root of $\mathfrak{F}$.

\item 
  \begin{itemize}
      \item For any  $w_i \in W_i$,
      $\text{end}_{\mathfrak{F}}(w_i)=\text{end}_{\mathfrak{F}_i}(w_i)$.

      \item $W_i \cap \text{End}(\mathfrak{F})=\text{End}(\mathfrak{F}_i)$.

      \item $ \text{End}(\mathfrak{F})=\text{End}(\mathfrak{F}_1) \cup \text{End}(\mathfrak{F}_2)$.

      \item For any  $ \langle w_1, w_2 \rangle \in W_1 \times W_2$, $\text{end}_{\mathfrak{F}}(\langle w_1, w_2\rangle )=\text{end}_{\mathfrak{F}_1}(w_1) \cup \text{end}_{\mathfrak{F}_2}(w_2)$.
  \end{itemize}

\item For any  $2 \leq n \in \omega$, any $u, x_0, \ldots, x_{n-1} \in W_i$, 

$u$ strongly unites $x_0, \ldots, x_{n-1}$ in $\mathfrak{F}_i \Leftrightarrow u$ strongly unites $x_0, \ldots, x_{n-1}$ in $\mathfrak{F}$.
\end{enumerate}
\end{lemma}

\begin{proof}
~
\begin{enumerate}

\item This follows directly from the definition of $R$. 
The third sentence can be proved using the second one.

\item This is straightforward, using the reflexivity and transitivity of $R_1, R_2$.

\item This is straightforward.

\item Part (1): This follows directly from part (1) of item 1.

Part (2): This follows directly from part (1) of item 1.

Part (3): By the reflexivity of $R_i$, it is easy to prove that $(W_1 \times W_2) \cap \text{End}(\mathfrak{F})=\emptyset$, and combining part (2), the goal is easily achieved.

Part (4): $\supseteq$: This follows from the reflexivity and transitivity of $R$ and part (2); $\subseteq$: This follows from parts (2) and (3) and the definition of $R$.

\item \quad Left

$\Leftrightarrow $ for any  $ i \in n$, $(uR_k x_i \text{\ and\ } u \in \Box_k(\blackdiamond_k\{x_i\} \cup \bigcup_{j \in n\{i\}}\Diamond_k \blackdiamond_k\{x_j\})$) \hfil\quad (definition of strong union)

$\Leftrightarrow$ for any  $i \in n$, $(uR x_i \text{\ and\ } u \in \Box(\blackdiamond\{x_i\} \cup \bigcup_{j \in n\{i\}}\Diamond \blackdiamond\{x_j\})$) \quad ($u, x_0, \ldots,$ $x_{n-1}$ $\in W_i$,  1)

$\Leftrightarrow$ Right \quad (definition of strong union).
\end{enumerate}
\qed  \end{proof}



\begin{lemma}\label{2. $u$ is a Hereditary Union Function}
The identity function on $W_1 \times W_2$ is a normal hereditary union function in $\mathfrak{F}_1 \otimes \mathfrak{F}_2$.
\end{lemma}

\begin{proof}
This is straightforward.
\qed  
\end{proof}

\begin{lemma}[connected product and $(su_n)$]\label{HU-merging Frame and (su_n)}
For any  $n \in \omega^*$,
if for any  $m\leq n$, $\mathfrak{F}_i$ satisfies $(su_m)$ ($i=1,2$), then $\mathfrak{F}_1 \otimes \mathfrak{F}_2$ satisfies $(su_n)$.
\end{lemma}

\begin{proof}
Let $\mathfrak{F}=\mathfrak{F}_1 \otimes \mathfrak{F}_2 =\langle W,R \rangle$. Assume the premise. We need to prove that $\mathfrak{F}$ satisfies $(su_{n})$.

We consider the following cases:
\begin{itemize}

\item For any  $k \in \{1,2\}$, any $w_k \in W_k$, and any $(w_kR)x_0,\ldots,x_{n-1}$,

by Lemma \ref{0. HU-merging Frame Basic Properties}.1, we have $w_k R_k x_0,\ldots,x_{n-1}$,

and since $\mathfrak{F}_k$ satisfies $(su_n)$, there exists  $(w_kR_k)z$, $z$ strongly unites $x_0,\ldots,x_{n-1}$ in $\mathfrak{F}_k$,

and by Lemma \ref{0. HU-merging Frame Basic Properties}.1 and Lemma \ref{0. HU-merging Frame Basic Properties}.5, there exists  $(w_kR)z$, $z$ strongly unites $x_0,\ldots,x_{n-1}$ in $\mathfrak{F}$.

\item For any  $w_1 \in W_1$, any $w_2 \in W_2$, and any $(\langle w_1, w_2\rangle R)x_0,\ldots,x_{n-1}$, we consider the following subcases:

\begin{itemize}

\item If there exists  $k \in \{1,2\}$, $x_0,\ldots,x_{n-1} \in W_k$,

then by Lemma \ref{0. HU-merging Frame Basic Properties}.1, $w_k R_k x_0,\ldots,x_{n-1}$, and since $\mathfrak{F}_k$ satisfies $(su_n)$, there exists  $(w_kR_k)z$, $z$ strongly unites $x_0,\ldots,x_{n-1}$ in $\mathfrak{F}_k$,

then by the definition of $R$ and Lemma \ref{0. HU-merging Frame Basic Properties}.5, $\langle w_1, w_2\rangle Rz$ and $z$ s.u. $x_0,\ldots,x_{n-1}$ in $\mathfrak{F}$.

\item Otherwise, there exist $i,j,l \in \omega$, such that $i+l \geq 1$ and $j+l \geq 1$ and $i+j+l=n$ and $x_0,\ldots,x_{n-1}$ can be arranged as ``$a_0,\ldots,a_{i-1}, \langle x_0, y_0\rangle ,\ldots$, $\langle x_{l-1},y_{l-1}\rangle , b_0,\ldots,b_{j-1}$'', where $a_0,\ldots,a_{i-1}, x_0,\ldots,x_{l-1} \in W_1$ and $b_0$, $\ldots,b_{n-1}, y_0,\ldots,y_{l-1} \in W_2$.

By $\langle w_1, w_2\rangle  R x_0,\ldots,x_{n-1}$ and Lemma \ref{0. HU-merging Frame Basic Properties}.1, we have $w_1 R_1 a_0,\ldots,a_{i-1}$, $x_0,\ldots,x_{l-1}$ and 
$w_2 R_2 b_0,\ldots,b_{n-1}, y_0,\ldots,y_{l-1}$.

Combining  Lemma \ref{0. HU-merging Frame Basic Properties}.5 and the fact that for any  $m\leq n, \mathfrak{F}_k$ satisfies $(su_m)$, $i,j \leq n$, we know that there exists  $(w_1R_1)u_1$, s.t. $u_1$ strongly unites $a_0,\ldots,a_{i-1}$, $x_0,\ldots,x_{l-1}$ in $\mathfrak{F}$, and there exists  $(w_2R_2)u_2$, s.t. $ u_2$ strongly unites $b_0,\ldots,b_{n-1}$, $y_0,\ldots,y_{l-1}$ in $\mathfrak{F}$.

Since the identity function on $W_1 \times W_2$ is a normal hereditary union function in $\mathfrak{F}$ (Lemma \ref{2. $u$ is a Hereditary Union Function}) and $R$ is reflexive and transitive (Lemma \ref{0. HU-merging Frame Basic Properties}.2), by Lemma \ref{-1. Hereditary Union Function and Strong Union} we have that $\langle u_1,u_2\rangle $ strongly unites $a_0,\ldots,a_{i-1}$, $\langle x_0, y_0\rangle ,\ldots,\langle x_{l-1},y_{l-1}\rangle , b_0,\ldots,b_{j-1}$ in $\mathfrak{F}$.

Since ``$a_0,\ldots,a_{i-1}, \langle x_0, y_0\rangle ,\ldots,\langle x_{l-1},y_{l-1}\rangle , b_0,\ldots,b_{j-1}$'' is a permutation of $x_0,\ldots,x_{n-1}$, it is easy to prove that $\langle u_1,u_2\rangle $ strongly unites $x_0,\ldots,x_{n-1}$ in $\mathfrak{F}$.

\end{itemize}
\end{itemize}
\qed  
\end{proof}

\begin{theorem}[Disjunction property]
For any  $\alpha,\beta\in Form$,
\[\vDash_{\textbf{SU}} \alpha_i\lor\alpha_2 \ \Rightarrow\ \vDash_{\textbf{SU}} \alpha_1\text{ or }\vDash_{\textbf{SU}} \alpha_2\]
\end{theorem}

\begin{proof}
By Lemma \ref{HU-merging Frame and (su_n)}, it is easy to prove that for any  rooted frames $\mathfrak{F}_1,r_1$ and $\mathfrak{F}_2,r_2$ in $\mathcal{C}_{su}$, there exists  a rooted frame $\mathfrak{F},r$ in $\mathcal{C}_{su}$ such that $\mathfrak{F}_i,r_i$ is isomorphic to two disjoint generated subframes of $\mathfrak{F}$. Thus, it is easy to prove that $\vDash_{\mathcal{C}_{su}}$, i.e., $\vDash_{\textbf{SU}}$, has the disjunction property.
\qed  
\end{proof}

~\\
~\\

%
%
%
\bibliographystyle{splncs04}
\bibliography{thesis}
\end{document}